\def\N{{\bf N}}
\def\Z{{\bf Z}}
\def\e{A}
\def\W{W}
\def\P{{\bf P}}
\def\A{{\bf A}}
\def\D{{\bf D}}
\def\C{{\bf C}}
\def\B{{\bf B}}

\input diagrams.tex 

\noindent
\bigskip

\centerline{\tenbf Symbolic dynamics and the category of graphs}
\vskip .25 true cm 

\smallskip
\vskip .25 true cm
\centerline{Terrence Bisson \quad\& \quad Aristide Tsemo}
\centerline{bisson@canisius.edu \quad \  \quad tsemo58@yahoo.ca}
\vskip .25 true cm

\bigskip\noindent
{\bf Abstract:} 
Symbolic dynamics is partly the study of walks in a directed graph. By a
walk, here we mean a  morphism to the graph from the Cayley graph of the monoid
of non-negative integers. Sets of these walks are also
important in other areas, such as stochastic processes, 
automata, combinatorial group theory, $C^*$-algebras, etc.
We put a Quillen model structure on the category of directed graphs,
for which the weak equivalences are those graph morphisms which induce
bijections on the set of walks.
We determine the resulting homotopy category.
We also introduce a ``finite-level'' homotopy category which
respects the natural topology on the set of walks.
To each graph we associate a basal graph, well defined up to isomorphism. We
show that the basal graph is a homotopy invariant for our model structure, and
that it is a finer invariant than the zeta series of a finite graph. We also
show that, for finite walkable graphs, if $B$ is basal and separated then the
walk spaces for $X$ and $B$ are topologically conjugate   if and only if  $X$
and $B$ are homotopically equivalent for our model structure.

\bigskip\noindent
{\bf 0. Introduction.}  
Symbolic dynamics is partly the study of walks in a directed graph; see the discussion in
Kitchens [1998] or Lind and Marcus [1995], for instance.
Sets of these walks are also
important in other areas, such as stochastic processes, 
automata, combinatorial group theory, $C^*$-algebras, etc.,
as can be seen from references such as Kemeny-Snell-Knapp [1976], Sakarovitch [2009], 
Epstein [1992], and Raeburn [2005] .

Let Gph denote the category of directed graphs.   In this paper we  investigate Gph
as a framework for analyzing symbolic dynamics of walks.
By a walk in directed graph $X$
we mean a  morphism from $N$ to $X$, where $N$
has a node $n$,
and  an arc from $n$ to $n+1$, for each natural number $n$.
So $N$ is a Cayley graph, of the following simple type. 
Any monoid $G$, together with some $a\in G$, generates a Cayley graph 
with a node for each element of $G$
and with an arc from $x$ to $xg$ for each node $x$.
Our results suggest that more general investigations of 
categories of $G$-sets and Cayley graphs are also useful,
but we leave that for further work.

\smallskip
In Section 1 we give our precise definitions and background.

\smallskip
In Section 2 we discuss the notion of Quillen model structure on a category,
which expedites the description of an associated homotopy category.
We define a model structure on Gph,
for which the weak equivalences are those graph morphisms which induce
bijections on the set of walks.

\smallskip
In Section 3
we determine the resulting homotopy category.

\smallskip
In Section 4 we describe the natural topology on the set of walks,
and  introduce a ``finite-level'' homotopy category of graphs which
respects the topology.

\smallskip
In Section 5 we explore some applications of covering morphisms,
inspired by the paper of Boldi and Vigna [2002].
We say that a graph is basal if the only epic covering morphisms defined on it
are the isomorphisms.
To each graph we associate a basal graph, well defined up to isomorphism.
We show that the basal graph is a homotopy invariant for our model structure,
and that it is a finer invariant than the zeta series of a finite graph.
We also show that, for finite walkable graphs,
if $B$ is basal and separated then the walk spaces 
for $X$ and $B$ are topologically conjugate   if and only if  $X$ and $B$ are homotopically
equivalent for our model structure.

\bigskip
The Quillen model on graphs that we investigate here seems
to be a particular example of the following general construction.
Let ${\cal E}$ be a topos, and ${\cal I}$ a family of objects of ${\cal E}$. 
A closed model can be defined on ${\cal E}$ for which the class of weak
equivalences are morphisms $f:X\rightarrow Y$ such that
${\rm Hom}_{\cal E}(i,X)\rightarrow {\rm Hom}_{\cal E}(i,Y)$ is a bijection for every $i\in {\cal I}$. 
In this paper, we study the particular example of this situation
when ${\cal E}$ is the topos of directed graphs and ${\cal I}$ has the single object $N$.
It seems likely that the general construction can be applied in 
other categories of combinatorial interest.

\medskip\noindent
{\bf 1. The set of walks and N-equivalence of graphs.} 

In this paper we continue our study of the category Gph of directed 
and possibly infinite graphs, 
with loops and multiple arcs allowed. 
This is the category studied in Bisson and Tsemo [2008], [2009].

Let us make precise the objects and morphisms in the category Gph. 
A  {\it graph} is a data-structure $X=(X_0,X_1,s,t)$ 
with a set $X_0$ of {\it nodes},
a set $X_1$ of {\it arcs}, and a pair of functions $s,t:X_1\to X_0$ which
specify the {\it source} and {\it target} node of each arc.
We may say that $a\in X_1$ is an arc
{\it from node $s(a)$ to node $t(a)$};
 a {\it loop} is just an arc $a$ with $s(a)=t(a)$.
A {\it graph morphism} $f:X\to Y$ is a pair of functions 
$f_1:X_1\to Y_1$ and $f_0:X_0\to Y_0$ such that
$s\circ f_1=f_0\circ s$ and $t\circ f_1=f_0\circ t$.
For $X$ and $Y$ in Gph, we may
sometimes denote the set of graph morphisms from $X$ to $Y$ by $[X,Y]$. 

The category Gph is very nice category to work with.  
In particular, it is a presheaf topos
(see Mac Lane and Moerdijk [1994], for instance, for a nice survey). 
As such, it has all limits and colimits, including the initial graph $0$ 
(with no nodes and no arcs)
and the terminal graph $1$ (with one node and one loop).  
Here are some other standard graphs that we will be using.
Let $\N$ denote the graph with nodes the natural numbers 
and arcs the pairs $(n,n+1)$ for $n\geq 0$,
with $s(n,n+1)=n$ and $t(n,n+1)=n+1$.  Let $\Z$ have nodes 
the integers and arcs $(n,n+1)$ for all integers,
with source and target as above.  Similarly,
let ${\bf P}_n$ have nodes $k$ and arcs $(k,k+1)$, for $0\leq k\leq n$.
We may call ${\bf P}_n$ the {\it path with $n$ arcs}, 
and use the notations ${\bf D}=\P_0$ and ${\bf A}=\P_1$.
For $n>0$, let ${\bf C}_n$ have the nodes the integers mod $n$, and arcs $(k,k+1)$, 
for $0\leq k\leq n$. 
We may call ${\bf C}_n$ the cyclic graph with $n$ arcs. Note that ${\bf C}_1=1$.

For any graph $X$, a {\it path of length $n$} is just a graph morphism $\alpha:{\bf P}_n\to X$;
its {\it source} $s(\alpha)$ is the image in $X$ of node $0$ in ${\bf P}_n$;
its {\it target} $t(\alpha)$ is the image in $X$ of node $n$ in ${\bf P}_n$.
Let $\alpha\beta$ denote the {\it concatenation} of paths, defined when $t(\alpha)=s(\beta)$.
We may denote the set of paths $[\P_n,X]$ by $P_n(X)$.

A {\it walk} $\omega$ in a graph $X$ 
is just a graph morphism $\omega:\N\to X$;
its {\it source} $s_0(\omega)$
is the image in $X$ of the node $0$ in $\N$.
Let $N(X)=[\N,X]$ denote the set of walks in $X$.
A graph morphism $f:X\to Y$ gives a natural function $N(f):N(X)\to N(Y)$ 
by $\omega\mapsto f\circ\omega$,
so that we have a functor from Gph to Set.

But in fact $N(X)$ is endowed with a natural {\it shift operation}, as follows.
Let $\tau:\N\to\N$ denote the graph morphism given on nodes by $\tau(n)=n+1$.
Let the shift operation $\tau: N(X)\to N(X)$  (with slight abuse of notation)  
be given by $\omega\mapsto\omega\circ\tau$ for $\omega\in N(X)$:
the shift of a walk just deletes the first arc in the walk. 
For any graph morphism $f$ the function $N(f)$  preserves $\tau$, in that
$N(f)\circ\tau=\tau\circ N(f)$.  So $N(X)$ is naturally an {\it N-set}, 
and we have a functor from Gph to NSet, in the following sense.

\smallskip\noindent
{\bf Definition}: An {\it N-set} is a pair $(S,\tau)$ 
with $\tau$ a function from $S$ to $S$; and
a {\it map of N-sets} from $(S,\tau)$ to $(S',\tau')$ is a function $f:S\to S'$ such that
$\tau'\circ f=f\circ\tau$.  
Let NSet denote the category of N-sets, with functor $N:{\rm Gph}\to{\rm NSet}$.
An  {\it N-equivalence} is a graph morphism $f:X\to Y$ for which
$N(f):N(X)\to N(Y)$ is an isomorphism of N-sets. 

\smallskip
There is a more general point of view about the category NSet.
Let $G$ be a monoid, with associative binary operation $G\times G\to G:(g,h)\mapsto g*h$ 
and with neutral element $e$; 
a $G$-set is a set $S$ together with an {\it action} that is a function  $\mu:G\times S\to S$ 
such that $\mu(e,x)=x$ and $\mu(g,\mu(h,x))=\mu(g*h,x)$.
For any monoid $G$, the category of $G$-sets is a presheaf category,
and thus a topos; see Mac Lane and Moerdijk [1994], for instance.
Then  NSet can be viewed as the category of actions of the
monoid $N$ of natural numbers, under addition, since a set $S$ together with an arbitrary 
function $\tau:S\to S$ corresponds
exactly to an action of the monoid $N$, 
by $\mu(n,x)=\tau^n(x)$ for $n\in N$.
Thus we can view NSet as a presheaf topos, with 
all products, and all coproducts (sums) formed ``elementwise'', etc.

\smallskip\noindent
{\bf Definition:} The {\it arc graph} $\e(X)$ of a graph $X$ is
the graph with the arcs of $X$ as its nodes,
and with length $2$ paths in $X$ as its arcs; and with source and target given by
$s(a_1,a_2)=a_1$ and $t(a_1,a_2)=a_2$.
Let $s_{1,0}:\e(X)\to X$ denote the graph morphism given on nodes by $a\mapsto s(a)$, and
on arcs by $(a',a)\mapsto a'$.
This is a graph morphism since
each arc $(a',a)$, from  $a'$ to $a$, in $\e(X)$ maps to the arc $a'$, from $s(a')$ to $s(a)=t(a')$, in $X$.

\smallskip
We will prove that $s_{1,0}:\e(X)\to X$ is an N-equivalence, in Section 4,
as part of a more general analysis..
The arc graph is sometimes called  ``the line digraph''  or ``the line graph for directed graphs'';
see for instance  Kotani and  Sunada [2000], where it i used in
connection with zeta series.

Here are some examples.  We have $\e(\P_n)=\P_{n-1}$; in particular, $\e(\D)=0$ and $\e(\A)=\D$.
Also, $\e(\N)=\N$ and $\e(\Z)=\Z$, and $\e(\C_n)=\C_n$; in particular, $\e(1)=1$.
For any set $S$, let ${\bf B}(S)$ denote the ``bouquet of loops'' with one node and with $S$
as its set of arcs.
Then $\e({\bf B}(S))={\bf K}(S)$ is the ``very complete graph'' with nodes $S$ and arcs $S^2$,
and with exactly one arc between any two nodes (including a unique loop from each node to itself).
The equal signs above are really denoting natural isomorphisms, of course.

Not every graph arises as an arc graph; for instance, $\e(X)$ is always a graph with no parallel arcs
(where two arcs $a$ and $a'$ with $s(a)=s(a')$ and $t(a)=t(a')$ are said to be {\it parallel}).

\medskip\noindent
{\bf 2. A model structure for N-equivalence of graphs.} 

In two previous papers (Bisson, Tsemo [2008], [2009])
we developed a Quillen model structure on the
category Gph, based on the set of cycles in a graph; we may refer to this
as the $C_*$-equivalence model, since here we 
will develop a different (simpler) Quillen model structure for Gph,
based on the set of walks in a graph.

We will use the following convenient terminology to explain Quillen model structures.
Let $\ell:X\to Y$ and $r:A\to B$ be morphisms in a category ${\cal E}$.
We say that $\ell$ is {\it weak orthogonal} to $r$ (abbreviated by $\ell\dagger r$) 
when all squares with $r$ on the right and $\ell$ on the left can be filled:
$${\rm if}\quad \diagram
          X & \rTo^{f} & A  \cr
   \dTo^{\ell} &           &  \dTo^{r}   \cr
          Y  & \rTo^g   &B
                                    \enddiagram \quad{\rm commutes,\ then}\quad
           \diagram[nohug]
          X & \rTo^{f} & A  \cr
   \dTo^{\ell} &  \NE^{h} &  \dTo^{r}   \cr
          Y & \rTo^g   &B
                                    \enddiagram\quad {\rm commutes\ for\ some\ } h.$$                                  
Given a class ${\cal F}$ of morphisms we define
${\cal F}^\dagger=\{r: f\dagger r, \ \forall f\in{\cal F}\}\quad{\rm and}\quad
{}^\dagger{\cal F}=\{\ell: \ell\dagger f, \ \forall f\in{\cal F}\}$.
A {\it weak factorization system} in ${\cal E}$ is given by two classes
${\cal L}$ and ${\cal R}$, such that ${\cal L}^\dagger={\cal R}$ and
${\cal L}={}^\dagger{\cal R}$ and such that, for any morphism $c$ in ${\cal E}$,
there exist $\ell\in{\cal L}$ and $r\in{\cal R}$ with $c=r\circ\ell$.

\smallskip 
We may express Quillen's notion [1967] 
of ``model category structure'' via the following axioms, 
which we learned from Section 7 of Joyal and Tierney [2007].

\smallskip\noindent
{\bf Definition:} A {\it model structure} on 
a category ${\cal E}$ with finite limits and colimits
is a triple $({\cal C},{\cal W},{\cal F})$
of classes of morphisms in ${\cal E}$ which satisfy

1) ``three for two'': if two of the three morphisms  $a, b, a\circ b$ belong to ${\cal W}$
then so does the third,

2) the pair $(\underline{\cal C} , {\cal F})$ is a weak factorization system 
(where $\underline{\cal C}={\cal C}\cap{\cal W}$),

3) the pair $({\cal C} , \underline{\cal F})$ is a weak factorization system 
(where $\underline{\cal F}={\cal W}\cap{\cal F}$).

\noindent
For instance, the {\it trivial} model structure (for any suitable category ${\cal E}$)
is given by the triple $({\rm All},{\rm Iso},{\rm All})$.

\noindent
The morphisms in  ${\cal W}$ are called {\it weak equivalences}.
The morphisms in  ${\cal C}$ are called {\it cofibrations}, and
the morphisms in  $\underline{\cal C}$ are called {\it acyclic cofibrations}.
The morphisms in  ${\cal F}$ are called {\it fibrations}, and
the morphisms in  $\underline{\cal F}$ are called {\it acyclic fibrations}.
An object $X$ in ${\cal E}$ is called {\it cofibrant} 
when $0\to X$ is in ${\cal C}$ (a cofibration), where $0$ is an initial object.
Dually,  $X$ is called {\it fibrant} when $X\to 1$ is in ${\cal F}$ (a fibration),
where $1$ is a terminal object.

\smallskip 
We will show that the following three morphism classes
give a model structure on the category Gph:

\item{$\bullet$} the fibrations are ${\cal F}_N={\rm All}$, the collection of all graph morphisms,

\item{$\bullet$} the weak equivalences are ${\cal W}_N$, the collection of all N-equivalences, and

\item{$\bullet$} the cofibrations are ${\cal C}_N={}^\dagger{\cal W}_N$.

\smallskip
In Appendix A we give a direct proof, using a ``small object'' argument,
that $({\cal C}_N,{\cal W}_N,{\cal F}_N)$  
is a model structure on Gph.  We may call it the {\it N-equivalence model structure on Gph};
the subscripts here are optional, but serve to distinguish these classes from 
the $C_*$-equivalence model structure from Bisson and Tsemo [2008], [2009]).

\smallskip
In this section we will show that $({\cal C}_N,{\cal W}_N,{\cal F}_N)$
is a model structure, by identifying it with a ``transport''
of the trivial model structure from the category NSet.
This will also show that the N-equivalence model structure is {\it cofibrantly generated}.
The transport will be along an {\it adjunction} (pair of adjoint functors) 
between Gph and NSet;
see Mac Lane [1971] for general background on when a pair of functors form an adjunction.
Section 2.1 in Hovey [1999], for example, has a nice discussion of cofibrant generation,
and other concepts which will be used in the following, such as
transfinite compositions, small object arguments,  etc.

Let ${\cal E}$ be a category with all limits and colimits.
Briefly, a model structure $({\cal C},{\cal W},{\cal F})$ 
on ${\cal E}$ is {\it cofibrantly generated}
when there are sets $I$ and $J$ of morphisms which generate ${\cal C}$ and $\underline{\cal C}$,
in the sense that
${}^\dagger(I^\dagger)={\cal C}$ and ${}^\dagger(J^\dagger)=\underline{\cal C}$;
thus we also have $I^\dagger=\underline{\cal F}$ and $J^\dagger={\cal F}$.
For a set $H$ of morphisms in ${\cal E}$, 
let ${\rm cell}(H)$ denote  the class of all transfinite compositions of pushouts of morphisms in $H$;
the morphims in ${\rm cell}(H)$ are called {\it relative $H$-cell complexes}.
For background and references on the proof of the following general result, see
Berger and Moerdijk [2003], for instance.

\smallskip\noindent
{\bf Transport Theorem:} Let ${\cal E}$ be a model category which is cofibrantly generated,
with cofibrations generated by  $I$ and  acyclic cofibrations generated by $J$.
Let ${\cal E}'$ be a category with all limits and colimits, 
and suppose that we have an adjunction
$$L:{\cal E}\rightleftharpoons{\cal E}':R\quad{\rm with}\quad
R({\rm cell}\ L(J))\subseteq {\cal W}.$$
Also, assume that the sets $L(I)$ and $L(J)$ each permit the small object argument.
Then there is a cofibrantly generated model structure on ${\cal E}'$
with generating cofibrations $L(I)$ and generating acyclic cofibrations $L(J)$.
Moreover, the model structure $({\cal C}',{\cal W}',{\cal F}')$  satisfies
$f\in{\cal W}'$ iff $R(f)\in{\cal W}$,
and $f\in{\cal F}'$ iff $R(f)\in{\cal F}$.

\smallskip
We apply the transport theorem with ${\cal E}$ as the category NSet, and with
${\cal E}'$ as the category Gph.  We use an adjunction
$$D:\hbox{NSet}\rightleftharpoons \hbox{Gph}:N$$
which plays a central role throughout this paper.
We have already defined the functor $N$.
For any N-set $(S,\tau)$, let $X=D(S,\tau)$ denote the graph with nodes $X_0=S$ and arcs $X_1=S$, 
where the source and target functions $s,t:X_1\to X_0$ are given by 
$s(x)=x$  and $t(x)=\tau(x)$ for each $x\in S$. 
Thus the elements in the N-set $S$ give the nodes and the arcs in the graph $X$,
and each arc $x$ has target $\tau(x)$ and source $x$;
we think of $\tau(x)$ as telling the unique ``target'' of each element $x$ in the N-set $S$.

It is easy to check directly that
$(D,N)$ is an adjoint pair of functors; the adjunction is also proved
in Bisson and Tsemo [2009], but
there we used the functor from NSet to Gph which assigned to $(S,\tau)$ the graph 
directed opposite to $D(S,\tau)$. 
Here we are directing our arcs in the way that seems natural
in graphical representation of dynamical systems 
(see Article III in Lawvere and Schanuel [1997], for instance).

\smallskip\noindent
{\bf Proposition:} The trivial model structure on NSet, 
when transported along the adjunction $(L,D)$,  gives
the N-equivalence model structure $({\cal C}_N,{\cal W}_N,{\cal F}_N)$ on Gph.
This model structure is cofibrantly generated by 
${\bf I}=\{{\bf i},{\bf j} \}$ and by  ${\bf J}=\{ {\bf 0} \}$,
where ${\bf i}:0\to \N$ and ${\bf j}: \N+\N\to \N$ are the initial and co-diagonal graph morphisms,
and ${\bf 0}$ is the identity graph morphism ${\bf 0}:0\to 0$.

\smallskip\noindent
{\bf Proof:} First we make precise our terminology for morphisms ${\bf i}$
and ${\bf j}$.  Any object $X$ in a category with coproducts has
{\it initial} morphism $0\to X$ (where  $0$ is the initial object), and 
{\it co-diagonal} morphism $X+X\to X$ (the morphism from the coproduct $X+X$ determined
by the pair of identity morphisms).  The category of N-sets has coproducts;
the initial object $0$ is the empty set.  We (temporarily) let ${\rm N}$ denote
the N-set of natural numbers with shift map $\tau(n)=n+1$, and consider
the sets  $I=\{{\rm i},{\rm j} \}$ and $J=\{ {\rm 0} \}$ of N-set maps, 
with initial N-set maps ${\rm 0}:0\to 0$ and ${\rm i}:0\to {\rm N}$, 
and co-diagonal N-set map ${\rm j}: {\rm N}+{\rm N}\to {\rm N}$.
We have   $J^\dagger=\hbox{All}$, so that ${}^\dagger(J^\dagger)=\hbox{Iso}$;
and we have $I^\dagger=\hbox{Iso}$, so that ${}^\dagger(I^\dagger)=\hbox{All}$.
This shows that the trivial model structure on NSet is cofibrantly generated.
The smallness conditions  in the Transport Theorem are automatically satisfied in our
presheaf categories (see the proof at Example 2.1.5 in Hovey [1999], for instance).
Now, let  ${\bf I}=D(I)$ and ${\bf J}=D(J)$; then ${\bf I}=\{{\bf i},{\bf j} \}$ and
${\bf J}=\{ {\bf 0} \}$.  So, every morphism in ${\rm cell}({\bf J})$ is a graph isomorphism,
and the Transport Theorem applies, since we have
$f\in {\rm cell}\ D(J)$ implies $N(f)\in{\cal W}$.
We immediately have ${\bf J}^\dagger=\hbox{All}={\cal F}_N$ and 
${}^\dagger({\bf J}^\dagger)=\hbox{Iso}=\underline{\cal C}_N$.
Moreover, the definitions (in terms of filling conditions) show that
${\bf I}^\dagger={\cal W}_N=\underline{\cal F}_N$, so that
${}^\dagger({\bf I}^\dagger)={}^\dagger{\cal W}_N={\cal C}_N$. 
It follows that our morphism classes $({\cal C}_N,{\cal W}_N,{\cal F}_N)$ are
cofibrantly generated by ${\bf I}$ and ${\bf J}$.  QED    

\smallskip\noindent
{\bf Definition:} A graph $X$ is a {\it dynamic graph} when every node in $X$ has
exactly one arc leaving it.  Let DGph denote the full subcategory of dynamic graphs.

Thus the dynamic graphs are those which are isomorphic
to $D(S,\tau)$ for some N-set $(S,\tau)$.

\smallskip\noindent
{\bf Proposition:} For the N-equivalence model structure on category Gph,
every graph morphism is a fibration, and every graph morphism
between dynamic graphs is a cofibration.
In particular, every graph is fibrant, and every dynamic graph is cofibrant.

\smallskip\noindent
{\bf Proof:}  As part of the definition of the N-equivalence model structure,
every graph morphism is a fibration.  We can use the transport definition
of the model structure to get partial information about the class cofibrations, as follows.
Let $I$ denote the set $\{{\rm i},{\rm j} \}$  of N-set maps,
as in the proof of the previous proposition.
We showed there that the cofibrations in our N-equivalence model are
generated by the set $D(I)$ of morphisms in Gph, so that 
${\rm cell}(D(I))\subseteq {\cal C}_N$.
Since the functor $D$ is a left adjoint, it preserves all colimits;
so $D({\rm cell}(I))\subset {\rm cell}(D(I))$. 
But  every map $f:S\to T$ of N-sets is in
${\rm cell}(I)$, as follows: let $S'=S+\sum_{x\in T}{\rm N}$; then $S\to S'$ is a pushout of 
a sum of copies of ${\rm i}$; and $S'\to T$ is a pushout of copies of ${\rm j}$
(this is just like the argument that 
all functions between sets are in ${\rm cell}(\{1+1\to 1,0\to 1\})$). 
It follows that $D({\rm cell}(I))$ is the class of graph morphisms between
dynamic graphs, and these are cofibrations.     QED

\smallskip
The adjunction $(D,N)$ assigns, to each graph morphism $D(S,\tau)\to X$,
an N-set map $(S,\tau)\to N(X)$ (called its {\it adjoint}).
The adjoint to the identity morphism $D(S,\tau)\to D(S,\tau)$ 
is a natural N-set map $(S,\tau)\to N(D(S,\tau))$,
which is  called the {\it unit} of the adjoint pair $(D,N)$, at the N-set $(S,\tau)$.
For every N-set, the unit $(S,\tau)\to N(D(S,\tau))$ is an isomorphism of N-sets,
since there is a unique walk starting at each node in a dynamic graph.
Note that an N-set map which  is a bijection is an N-set isomorphism.

The natural graph morphism $D(N(X))\to X$ which is adjoint to the identity N-set map $N(X)\to N(X)$ 
is called the {\it counit} of the adjoint pair $(D,N)$, at the graph $X$.
We may refer to $\W(X)=D(N(X))$ as the {\it walk graph} of $X$;
it is the dynamic graph which 
has the walks in $X$ as both its nodes and its arcs, with
$s(\omega)=\omega$ and $t(\omega)=\tau(\omega)$, for $\omega$ any walk in $X$.
Then the counit of the adjunction is the graph morphism $s_0:\W(X)\to X$ which,
on nodes, assigns to each walk $\omega$ its first node; and on arcs assigns
to $\omega$ its first arc.  
We may refer to $s_0$ as the {\it source truncation}.

\smallskip\noindent
{\bf Proposition:}  For any graph $X$, the graph $\W(X)$ 
is cofibrant and the graph morphism $s_0:\W(X)\to X$ is an N-equivalence.
Also, $\W(f):\W(Y)\to\W(X)$ is a graph isomorphism for any N-equivalence $f:Y\to X$.

\smallskip\noindent
{\bf Proof:} Since $\W(X)$ is a dynamic graph, it is cofibrant.
Also $\W(X)\to X$ is an N-equivalence, since 
$N(\W(X))=N(D(N(X))=N(X)$, through the identification $N(D(S,\tau))=(S,\tau)$ for every N-set $(S,\tau)$.
The second statement follows from the fact that $D(N(f))$ is an isomorphism when $N(f)$ is an isomorphism.
QED

\smallskip
The above proposition shows that $\W:\hbox{Gph}\to\hbox{Gph}$ 
is the {\it coreflection} of Gph into the full subcategory DGph. 
See Mac Lane [1971] for definitions of the general concepts. 
Results in Bisson and Tsemo [2009] show, essentially, that
DGph is a full reflective and coreflective subcategory of Gph.

\smallskip\noindent
{\bf Corollary:} The dynamic graphs are the cofibrant objects 
for the N-equivalence model structure on graphs.

\smallskip\noindent
{\bf Proof:} We have already shown that every dynamic graph is cofibrant.  For the converse,
suppose that graph $X$ is a cofibrant graph.  
Since $s_0:\W(X)\to X$ is an N-equivalence, we have a filling $f$ for the diagram
$$\diagram[nohug]
    0     & \rTo      &\W(X)     \cr
     \dTo & \NE^{f}   &\dTo^{s}  \cr
    X & \rTo^{\rm id} & X        \cr
                                       \enddiagram$$
This implies that $s$ is an epic graph morphism and that $f$ is a monic graph morphism.
Suppose that $X$ is not a dynamic graph;
then the set $X(x,*)$ of  arcs leaving some node $x$ in $X$ has
cardinality other than one. 
But $X(x,*)$ can't be empty, since then there would be no walk in $X$ leaving $x$,
and $x$ would not be in the image of $s_0:\W(X)\to X$, which contradicts $s$ being epic.
So $X(x,*)$ must have more than one element.
But $\W(X)$ is a dynamic graph, so
every arc in $X(x,*)$ must map to the unique arc leaving $f(x)$ in $\W(X)$, which contradicts
$f$ being monic.  
QED

\medskip\noindent
{\bf 3. The N-equivalence homotopy category.}

The purpose of giving a model structure on a category ${\cal E}$ is to construct 
and study a new category ${\rm Ho}({\cal E})$
which inverts the weak equivalences of the model category.  Let us explain. 

Suppose that ${\cal E}$ is a model category.
A functor with domain ${\cal E}$
is said to be a {\it homotopy functor} when it takes every $f\in{\cal W}$ to an isomorphism.
This involves just the class ${\cal W}$ of weak equivalences in the model structure.
Quillen [1967] used the classes ${\cal C}$ and ${\cal F}$ to
describe a  particular category ${\rm Ho}({\cal E})$,
together with a functor $\gamma:{\cal E}\to {\rm Ho}({\cal E})$ which is
{\it initial} among homotopy functors on ${\cal E}$.
This means that $\gamma$ is a homotopy functor and that
any homotopy functor $\Phi:{\cal E}\to{\cal D}$ factors 
uniquely through $\gamma$, in that $\Phi=\Phi'\circ \gamma$ for a unique functor 
$\Phi':{\rm Ho}({\cal E})\to {\cal D}$.  

In fact, Quillen constructs the category  ${\rm Ho}({\cal E})$ to have
the same objects as ${\cal E}$, and describes
the set ${\rm Ho}(X,Y)$ of ``homotopy arrows''  from $X$ to $Y$
in ${\rm Ho}({\cal E})$, for any objects $X$ and $Y$ in ${\cal E}$.
His construction uses the following notions.
A {\it cofibrant replacement} for an object $X$ in ${\cal E}$ is a 
morphism $f:X'\to X$ where $X'$ is cofibrant and $f$ is a weak equivalence and a fibration
($f\in\underline{\cal F}={\cal W}\cap{\cal F}$).
Dually, a {\it fibrant replacement} for $X$  is a 
morphism $g:X\to X''$ where $X''$ is fibrant and $g$ is a weak equivalence and a cofibration
($g\in\underline{\cal C}={\cal W}\cap{\cal C}$).
It follows from the model category axioms that each object in ${\cal E}$ has a
cofibrant replacement and a fibrant replacement.

Then the homotopy functor $\gamma:{\cal E}\to {\rm Ho}({\cal E})$ carries morphisms 
in ${\cal E}$ to homotopy arrows in ${\rm Ho}({\cal E})$,
but there are usally homotopy arrows in ${\rm Ho}({\cal E})$ which are not equal to 
$\gamma(f)$ for any morphism $f$ in ${\cal E}$.  
So morphisms in ${\cal E}$ may become
invertible in ${\rm Ho}({\cal E})$, and
objects which are not isomorphic in ${\cal E}$ may become
isomorphic in ${\rm Ho}({\cal E})$.
We may say that two objects $X$ and $Y$ in ${\cal E}$ are {\it homotopy-equivalent}
when $X$ and $Y$ become isomorphic in ${\rm Ho}({\cal E})$;
and that a morphism $f:X\to Y$ in ${\cal E}$ is a {\it homotopy equivalence}
when $\gamma(f)$ becomes invertible in ${\rm Ho}({\cal E})$. 
Also, we may say that morphisms $f,g:X\to Y$ in ${\cal E}$  are {\it homotopic}
when they become equal in ${\rm Ho}({\cal E})$, with $\gamma(f)=\gamma(g)$.
Quillen's description of the homotopy arrows ${\rm Ho}({\cal E})$ uses the following notions.

\smallskip
Let us see how these ideas work out for our N-equivalence model structure on Gph.
Recall that every graph morphism is a fibration and that every graph is fibrant;
every graph is its own fibrant replacement.  Moreover,
our results at the end of section 2 show that the natural graph morphism $s_0:\W(X)\to X$ 
gives a cofibrant replacement for every graph $X$.

\smallskip\noindent
{\bf Proposition:} The functor $N:\hbox{Gph}\to \hbox{NSet}$ 
induces an equivalence of categories $\hbox{Ho(Gph)}\to \hbox{NSet}$. 

\smallskip\noindent
{\bf Proof:} We show that $N:\hbox{Gph}\to \hbox{NSet}$  factors 
through $\gamma:\hbox{Gph}\to \hbox{Ho(Gph)}$.
The functor $N:\hbox{Gph}\to \hbox{NSet}$ factors through 
$\gamma:\hbox{Gph}\to \hbox{Ho(Gph)}$ and $N:\hbox{Ho(Gph)}\to \hbox{NSet}$,
which gives the desired equivalence.  Note that the unit
$N(D(S,\tau))\to (S,\tau)$ is already an isomorphism
and it is only necessary to recall that the N-equivalence $\W(X)\to X$
can be viewed as the counit $D(N(X))\to X$. QED

\smallskip
For any graph $X$, consider the subgraph of $X$ which is
the image of the natural graph morphism $s_0:\W(X)\to X$.
We will call it the {\it walkable subgraph} of $X$.
Now we are ready to describe precisely the various notions of homotopy for 
the N-equivalence model structure on Gph.

\smallskip\noindent
{\bf Proposition:}
Graphs $X$ and $Y$ are homotopy-equivalent
if and only if the N-sets $N(X)$ and $N(Y)$ are isomorphic.
A graph morphism $f$ is  a homotopy equivalence if and only if it is an N-equivalence.
Graph morphisms $f,g:X\to Y$ are homotopic 
if and only if they agree on the walkable subgraph of $X$.

\smallskip\noindent
{\bf Proof:}  The first statement follows from the previous proposition: objects $X$ and $Y$ are isomorphic in
${\rm Ho(Gph)}$ if and only if $N(X)$ and $N(Y)$ are isomorphic in NSet.
For the second statement, we use the following general result.
From Quillen's description of the category ${\rm Ho}({\cal E})$, 
for any model structure  $({\cal C},{\cal W},{\cal F})$,
it follows that $\gamma(f)$ is invertible in ${\rm Ho}({\cal E})$
if and only if $f$ is in ${\cal W}$ (see Hovey [1999], Theorem I.2.10, for instance).
So, a graph morphism $f:X\to Y$ has  $\gamma(f)$ invertible in ${\rm Ho(Gph)}$
if and only if $N(f)$ is an isomorphism of N-sets; 
and these N-equivalences are taken to form the class ${\cal W}_N$ 
of weak equivalences for our N-model structure on Gph.
So $f$ is an N-equivalence if and only if it is a homotopy equivalence.
Our proof of the third statement uses the following lemma.

\smallskip\noindent
{\bf Lemma:} The natural map $s_*:[\W(X),\W(Y)]\to[\W(X),Y]$ 
given by $s_*(f)=s\circ f$ is a bijection.

\smallskip\noindent
{\bf Proof of Lemma:} The adjoint pair $(D,N)$ gives a natural bijection
$${\rm NSet}[N(X),N(Y)]\cong[D(N(X)),Y].$$ 
We showed that the counit of the adjoint pair $(D,N)$ gives
a natural identification between $N\circ D$ and the identity functor;
it follows that the functor $D$ gives a natural bijection 
$${\rm NSet}[N(X),N(Y)]\cong[D(N(X)),D(N(Y))].$$
Recall that $\W=D\circ N$. The resulting bijection 
$[D(N(X)),D(N(Y))]\cong[D(N(X)),Y]$ can be identified with $s_*:[\W(X),\W(Y)]\to[\W(X),Y]$. QED

\smallskip\noindent
{\bf Proof of proposition, continued:} 
We have shown that $f$ and $g$ are homotopic if and only $N(f)=N(g)$.
The lemma shows that $N(f)=N(g)$ if and only if the graph morphisms
$s\circ\W(f),s\circ\W(g):\W(X)\to Y$ are equal.  
Let $s:\W(X)\to w(X)$ denote the epic graph morphism onto
image of the graph morphism  $s_0:\W(X)\to X$.  Then $s\circ\W(f)=f_|\circ s$,
where $f_|$ denotes $f$ restricted to $w(X)$.  So,  $s\circ\W(f)=s\circ\W(g)$ if and only if
$f_|\circ s=g_|\circ s$, which is equivalent to $f_|=g_|$ since $s$ is epic. QED

\smallskip
By the above, any graph is homotopy equivalent to its walkable subgraph.
So, if a graph $X$ has no walks, then $N(X)$ empty, and the walkable subgraph of $X$ is empty;
in this case, $X$ is homotopy equivalent to $0$,
and any two graph morphisms from $X$ to $Y$ are homotopic (for any graph $Y$).  
In particular, the graphs $O$ and $1$ are homotopy equivalent. For example, 
let $X$ have nodes $x,x_1,x_2$ with arcs $a_i$ from $x$ to $x_i$; let $Y$ have
nodes $y,y_1,y_2$ with arcs $b_i$ from $x_i$ to $x$.  
Then $X$ and $Y$ are homotopy equivalent even though there is no graph morphism between $X$ and $Y$.

\smallskip
A functor $F$ defined on Gph will be a homotopy functor for the
N-equivalence model structure if and only if $F(f):F(X)\to F(Y)$ is an isomorphism whenever
$f:X\to Y$ is an N-equivalence. 
For instance, the functor $\gamma:\hbox{Gph}\to\hbox{Ho(Gph)}$ is initial among homotopy functors;
and it is equivalent to the functor $N:\hbox{Gph}\to \hbox{NSet}$.
This also shows that the cofibrant replacement functor $\W:\hbox{Gph}\to \hbox{NSet}$ is a homotopy functor,
since $\W$ is $D\circ N$, and composing a homotopy functor with another functor gives a homotopy functor.

\smallskip\noindent
{\bf Proposition:} Let $F$ be a dynamic graph:

a) there is a natural graph morphism $\sigma:F\to F$ determined by $s(\sigma(a))=t(a)$ on arcs;

b) the functor from Gph to NSet given by $X\mapsto ([F,X],\sigma^*)$, with $\sigma^*(f)=f\circ\sigma$,
is a homotopy functor.

\smallskip\noindent
{\bf Proof:} 
We may identify $F=D(S,\tau)$ for some N-set $(S,\tau)$.
The function $\tau$ is in fact an N-set map $\tau:(S,\tau)\to (S,\tau)$,
and gives a graph morphism $D(\tau):D(S,\tau)\to D(S,\tau)$.  This gives $\sigma:F\to F$,
and a functor $F$ from Gph to NSet, with  $F(X)=([F,X],\sigma^*)$.  
We must show that if a graph morphism $f:X\to Y$ 
is an N-equivalence then $F(f)$ is an isomorphism of N-sets. 
But $F=D(S,\tau)$, and the adjunction $(D,N)$ shows that $F(X)$ can be identified with 
the set of N-set maps from $(S,\tau)$ to $N(X)$,
so that the functor $X\mapsto F(X)$ factors through $N:\hbox{Gph}\to \hbox{NSet}$. 
QED

\smallskip
For example, the functor $Z:\hbox{Gph}\to \hbox{NSet}$ given by $X\mapsto [\Z,X]$ is a homotopy functor,
since $\Z=D(Z,+1)$ is the dynamic graph with nodes the integers.  We may refer to elements of $[\Z,X]$ 
as two-way walks in $X$.

As another example, for any $n>0$ the functor $\hbox{Gph}\to \hbox{NSet}$ given by $X\mapsto [\C_n,X]$ 
is a homotopy functor, since ${\bf C}_n=D(Z/n,+1)$ is the dynamic graph with nodes the integers mod $n$.  
It follows that the functors $C_n:\hbox{Gph}\to \hbox{Set}$, with $C_n(X)=[\C_n,X]$,
are homotopy functors.
We refer to elements of $[\C_n,X]$  as cycles of length $n$ in $X$;
they can be identified with the set of $\omega\in N(X)$ such that $\tau^n(\omega)=\omega$.
For a {\it finite graph} $X$ (finitely many nodes and arcs),
the {\it zeta series} of $X$ is the formal power series
$${\rm Zeta}(u)={\rm exp}(\sum_{m=1}^\infty c_m{u^m\over m}),$$
where $c_m=|C_m(X)|$ for $m>0$.

\smallskip\noindent
{\bf Corollary:}  If $X$ and $Y$ are N-equivalent finite graphs 
then they have the same zeta series.

\smallskip
Let us say that a graph morphism $f:X\to Y$ is {\it acyclic}
when $C_n(f):C_n(X)\to C_n(Y)$ is a bijection for every $n>0$. 
In Bisson and Tsemo [2009], we studied the homotopy category of graphs
that results when one inverts the acyclic graph morphisms;
here we will call it the acyclic model structure on Gph.
Our main result in that paper
said that $X$ and $Y$  have the same zeta series if and only if
they are homotopy equivalent in the acyclic model structure.
Let us write $X\sim_{C} Y$ for this situation, and
write $X\sim_{N} Y$ when $X$ and $Y$ are homotopy equivalent for the N-equivalence 
model structure.

\smallskip\noindent
{\bf Proposition:} If $X$ and $Y$ are finite graphs, then  $X\sim_{N} Y$ implies $X\sim_{C} Y$.

\smallskip\noindent
{\bf Proof:} If  $X\sim_{N} Y$ then there is an isomorphism of N-sets $\phi:N(X)\to N(Y)$.
For each $n> 0$ this restricts to give a bijection $\phi:C_n(X)\to C_n(Y)$.
These are finite sets if $X$ and $Y$ are finite graphs; and then we have $c_n(X)=c_n(Y)$
for all $n> 0$.  Thus $X$ and $Y$ have the same zeta series, so that we have $X\sim_{C} Y$. QED

\smallskip 
In section 6 we give an example of finite graphs $X$ and $Y$ 
which have the same zeta function but are not N-equivalent, so that
we have $X\sim_{C} Y$ but not $X\sim_{N} Y$.

\smallskip 
Many other natural functors from Gph to Set are {\it not} homotopy functors.
For instance $X\mapsto [{\bf D},X]=X_0$ is not a homotopy functor, since $)$ and $1$ are homotopy 
equivalent graphs, but $[{\bf D},0]\neq[{\bf D},1]$.  Similar reasoning applies to 
$X\mapsto\pi_0(X)$, the set of components of the graph $X$,
formed as the coequalizer of the functions $s,t:X_1\to X_0$.

\medskip\noindent
{\bf 4. Arc graphs and finite-level homotopy.}

In section 1 we defined the arc graph  $\e(X)$ for any graph $X$.
In Section 3 we defined the walk graph $W(X)$ and showed that it
provides a cofibrant replacement for the N-equivalence model structure.
Here we extend and relate these constructions,
by the following general considerations.

Any pair of arrows $i_s,i_t:E_0\to E_1$ in a category ${\cal E}$
gives a {\it representable functor} $E_*:{\cal E}\to{\rm Gph}$, 
by assigning, to any object
$E\in {\cal E}$, the graph $E_*(X)$ with ${\cal E}[E_0,E]$ as set of nodes $E_0(X)$, 
and with ${\cal E}[E_1,E]$ as set of arcs $E_1(X)$; the source and target of
any arc $\alpha:E_1\to E$ in $E_1(X)$ are given by $s(\alpha)=\alpha\circ i_s$
and $t(\alpha)=\alpha\circ i_t$.  Here ${\cal E}[E',E]$ denotes the set of morphisms  from object
$E'$ to object $E$ in the category ${\cal E}$.

For instance, our cofibrant replacement functor $\W:{\rm Gph}\to {\rm Gph}$ comes in this way from
$i_s,i_t:{\bf N}\to {\bf N}$ in Gph, where $i_s$ is the identity graph morphism and $i_t$ is the shift
graph morphism. 

For each $n\geq 0$ we define a functor $\e^n:{\rm Gph}\to {\rm Gph}$ by
the pair $i_s,i_t:{\bf P}_n\to {\bf P}_{n+1}$, where
the graph morphisms are given on nodes by $i_s(k)=k$ and by $i_t(k)=k+1$.  We might 
refer to the graph morphisms $i_s$ and $i_t$ as the inclusion at the source of the path and at the target
of the path.   For $n=0,1$ we have natural isomorphisms $\e^0(X)=X$ and $\e^1(X)=\e(X)$,
from ${\bf D}={\bf P}_0$ and ${\bf A}={\bf P}_1$.

\smallskip
Returning to general considerations, if $i'_s,i'_t:E'_0\to E'_1$ in ${\cal E}$
is giving another representable graph functor, then any pair of arrows $f_0:E_0\to E'_0$
and $f_1:E_1\to E'_1$ such that $f_1\circ i_s=i'_s\circ f_0$ and $f_1\circ i_t=i'_t\circ f_0$
determines a representable natural transformation 
from functor $E_*$ to functor $E'_*$.

For instance, the natural graph morphism $s_0:\W(X)\to X$
comes from $f_0:{\bf P}_0\to {\bf N}$ and $f_1:{\bf P}_1\to {\bf N}$.
More generally, for each $n\geq 0$ we define natural graph morphisms $s_n:\W(X)\to \e^n(X)$ by
$f_0:{\bf P}_{n}\to {\bf N}$ and $f_1:{\bf P}_{n+1}\to {\bf N}$;
and for $n,m\geq 0$ we define natural graph morphisms
$s_{m,n}:\e^{n+m}(X)\to \e^n(X)$ by $f_0:{\bf P}_{n}\to {\bf P}_{n+m}$ 
and $f_1:{\bf P}_{n+1}\to {\bf P}_{n+m+1}$. In all these cases, the graph morphisms $f_i$ are
determined by the condition that they take node $0$ to node $0$.
We may call $s_n$ and $s_{m,n}$ the length $n$ ``source truncations''.
In particular, 
$s=s_0:\W(X)\to X$ and we have $s_{m,0}:\e^m(X)\to X$.

\smallskip\noindent
{\bf Proposition:}  For any graph $X$ we have

1)  $s_{m,n}\circ s_{n+m}=s_n:\W(X)\to \e^n(X)$ 
and $s_{m+k,n}=s_{m,n}\circ s_{k,n+m}:\e^{n+m+k}(X)\to \e^{n}(X)$.

2) $\W(X)={\rm lim}_n\ \e^n(X)$.

3) $\e^n(\e^m(X))=\e^{n+m}(X)$

4) $\e^n(\W(X))=\W(X)=\W(\e^n(X))$

\smallskip\noindent
{\bf Proof:} For part 1, we check compatibility of the representing graph morphisms.
For part 2, we verify the universal limit condition for the representing graph morphisms
$\P_n\to \N$.  For 3, we use that every path of length $n+m$ is uniquely the concatenation of
a path of length $n$ and a path of length $m$.  The following lemma
shows that the natural
graph morphisms $\W(s_{n,0}):\W(\e^n(X))\to \W(X)$ and 
$s_{n,0}:\e^n(\W(X))\to \W(X)$ are graph isomorphisms, proving part 4:

\smallskip\noindent
{\bf Lemma:} If $Y$ is a dynamic graph then $s_n:\W(Y)\to \e^n(Y)$ and $s_{m,n}:\e^{n+m}(Y)\to\e^n(Y)$
are graph isomorphisms.

\smallskip\noindent
{\bf Proof of Lemma:} We use the fact that 
a graph morphism between dynamic graphs is a graph isomorphism if and only if it is bijective on nodes.
This is true since any graph morphism between dynamic graphs has
the form $D(f)$ for some N-set map $f:S_1\to S_2$; but an N-set map is an isomorphism
in and only if it is a bijection on elements, and elements in $S$ correspond to nodes in $D(S)$.
Then we note that $s_0$ is a graph morphism between dynamic graphs;
and it is clearly bijective on nodes. The other parts are similar.  QED

\smallskip
By part 3 of the proposition,
we may think of $\e^n$ as an iterated composition of the functor $\e$ with itself,
and we may refer to $\e^n(X)$ as the {\it $n$-fold, or length $n$, arc graph} on $X$.
We also extend our examples of N-equivalences as follows.

\smallskip\noindent
{\bf Corollary:} The natural graph morphisms $s_n:\W(X)\to \e^n(X)$ and 
$s_{m,n}:\e^{m+n}(X)\to \e^n(X)$ are N-equivalences.

\smallskip\noindent
{\bf Proof:} We can see that $s_n:\W(X)\to \e^n(X)$ is an N-equivalence by identifying it with 
$\W(\e^n(X))\to \e^n(X)$ (using $\W(\e^n(X))=\W(X)$). Then $\e^n(X)\to X$ is an N-equivalence by 
the 2/3 property for N-equivalences (this could also be shown by induction on $n$, of course).
Finally, $\e^{m+n}(X)\to \e^n(X)$ is an N-equivalence, 
since $\W(X)\to \e^{m+n}(X)$ and $\W(X)\to \e^n(X)$
are N-equivalences. QED

\smallskip 
Recall that $\W(X)\to X$ gives a cofibrant replacement for our model structure,
and every graph is its own fibrant replacement.  We will show how homotopy arrows from $X$ to $Y$
are represented by graph morphisms from $\W(X)$ to $Y$.
In the following, recall that
$s_0:\W(\W(X))\to \W(X)$ is a graph isomorphism, for any graph $X$,
so that $\gamma(s_0):\W(X)\to X$ is an isomorphism in Ho(Gph).

\smallskip\noindent
{\bf Definition:} 
For graph morphisms $f:\W(X)\to Y$ and $g:\W(Y)\to Z$, let $g\odot f$ denote
the graph morphism $g\circ \W(f)\circ s_0^{-1}$:                        
$$\diagram
    \W(\W(X)) & \rTo^{\W(f)}  &\W(Y)     \cr
     \dTo^{s_0} &                &\dTo^{g}   \cr
          \W(X) & \rTo^{g\odot f}       & Z     
                                                    \enddiagram$$
For any graph morphism $f:\W(X)\to Y$, let $\gamma'(f)$ denote
$\gamma(f)\circ \gamma(s_0)^{-1}$ in Ho(Gph):  
$$ \diagram[nohug]
        \W(X) & \rTo^{\gamma(f)}        & Y  \cr
   \dTo^{\gamma(s_0)} &\NE_{\gamma'(f)} &    \cr
            X &                 &    \cr \enddiagram$$
\smallskip\noindent
{\bf Proposition:} 
The function  $\gamma'$ gives a bijection between the set of graph morphisms from $\W(X)$ to $Y$, 
and the set of homotopy arrows from $X$ to $Y$ in Ho(Gph).  
For graph morphisms $f:\W(X)\to Y$ and $g:\W(Y)\to Z$, 
we have $\gamma'(g\odot f)=\gamma'(g)\circ \gamma'(f)$ in Ho(Gph).

\smallskip\noindent
{\bf Proof:}
We use $[X,Y]$ as notation for the set of graph morphisms from $X$ to $Y$, etc.
In section 3 we showed the equivalence of Ho(Gph) and NSet, giving
natural bijections  
$$[\W(X),\W(Y)]\cong {\rm NSet}[N(X),N(Y)].$$
Let WGph denote the category with the same objects as Gph, 
but with the new set of 
morphisms $${\rm WGph}[X,Y]=[\W(X),\W(Y)]$$ for objects $X$ and $Y$.
The functor $\hbox{Gph}\to\hbox{WGph}$ given by $f\mapsto \W(f)$ is
a homotopy functor; in fact, $\W(f)$ is a graph isomorphism if and only if $N(f)$ is an isomorphism.
It follows that  Ho(Gph) and WGph are isomorphic as categories.
The natural bijection 
$$[\W(X),\W(Y)]\cong [\W(X),Y]$$
allows us also to describe Ho(Gph) as the category whose objects are the graphs, 
but with morphism sets ${\rm Ho}(X,Y)=[\W(X),Y]$.
Then the homotopy functor $\gamma:\hbox{Gph}\to \hbox{Ho(Gph)}$
is described by the natural functions 
$s_0^*:[X,Y]\to[\W(X),Y]$, where $s_0^*(f)= f\circ s_0$.
The composition in the  category Ho(Gph) corresponds to
as the associative ``composition'' 
$$(f,g)\mapsto g\odot f\quad\quad [\W(X),Y]\times [\W(Y),Z]\to [\W(X),Z].$$
QED

\smallskip  
Thus the category Ho(Gph) has been described
directly in terms of graph morphisms defined on dynamic graphs, which
are the cofibrant objects for our model structure.  For this reason, we
think of the above as giving a ``cofibrant description of the homotopy category''.
As an application, we note that a graph morphism $f:X\to Y$ is a homotopy equivalence 
if and only if there exists a graph morphism $q:\W(Y)\to X$  
(thought of as a ``homotopy arrow from $Y$ to $X$''),
with $f\odot q={\rm id}$ and $q\odot f={\rm id}$.
This says that $q$ makes the following diagram commute:  
$$\diagram[nohug]
          \W(X) & \rTo^{W(f)} & \W(Y)  \cr
   \dTo^{s} &  \SW^{q} &  \dTo^{s}   \cr
          X & \rTo^f   &Y
                                    \enddiagram$$ 
The existence of such a $q$ also shows that $\W(f):\W(X)\to \W(Y)$ is an isomorphism of graphs.  

\smallskip
The above cofibrant description of Ho(Gph) suggests the following
notion of ``finite-level homotopy''.

\smallskip\noindent
{\bf Definition.}
Define $\gamma_n:[\e^n(X),Y]\to[\W(X),Y]={\rm Ho}(X,Y)$ by $\gamma_n(f)= f\circ s_n$,
where $s_n:\e^n(X)\to X$.  A homotopy arrow from $X$ to $Y$ in Ho(Gph) is a
{\it homotopy arrow of  level $n$} when it has the form $\gamma_n(f)$ for some 
graph morphism $f:\e^n(X)\to Y$.
Letting $n$ vary gives the {\it finite-level homotopy arrows}.

\smallskip\noindent
{\bf Proposition:}
The finite-level homotopy arrows form a subcategory of Ho(Gph).

\smallskip\noindent
{\bf Proof:}
The identity graph morphisms  are homotopy arrows of level $0$, 
by the identification  $\e^0(X)= X$.
Consider the functions $[\e^n(X),Y]\times [\e^m(Y),Z]\to [\e^{m+n}(X),Z]$, 
defined by  $(f,g)\mapsto g\circ \e^m(f)$ for $f:\e^n(X)\to Y$ and $g:\e^m(Y)\to Z$ and
$\e^m(f):\e^{m+n}(X)\to \e^m(Y)$.  These give a ``composition'' which is
compatible with the composition in Ho(Gph), by the natural graph morphisms from walk graphs to arc graphs.
This shows that the finite-level homotopy arrows are closed under composition, and form a subcategory
of Ho(Gph).  QED

\smallskip
We may call this the {\it finite-level subcategory} of Ho(Gph).
Let us say that a graph morphism $f:X\to Y$ is 
a {\it level $n$ homotopy equivalence} when
there exists a graph morphism $q:\e^n(Y)\to X$  which fills the diagram: 
$$\diagram[nohug]
          \e^n(X) & \rTo^{\e^n(f)} & \e^n(Y)  \cr
   \dTo^{s} &  \SW^{q} &  \dTo^{s}   \cr
          X & \rTo^f   &Y
                                    \enddiagram$$
 If $f$ is a  level $n$ homotopy equivalence and $g$ is 
a level $m$ homotopy equivalence then $f\circ g$ is a 
level $n+m$ homotopy equivalence.  Also, if $f$ is a level $n$ homotopy equivalence then 
$f$ is a level $n+1$ homotopy equivalence.
For example, for every $n,m\geq 0$, the
graph morphism $s_{n,m}:\e^{n+m}(X)\to \e^m(X)$ is a level $n$ homotopy equivalence.
In particular, $s_{1,0}:\e(X)\to X$ is a level $1$ homotopy equivalence.

\smallskip  
Recall that any homotopy equivalence of graphs corresponds to an isomorphism of N-sets;
and we have picked out a subcategory of finite-level homotopy arrows and finite-level homotopy equivalences.
In the next section we show that a finite-level homotopy equivalence 
corresponds to a special kind of N-isomorphisn, called a ``topological conjugacy''.

\medskip\noindent
{\bf 5. Symbolic dynamics and topological conjugacy of walk spaces.} 

In this section we want to relate our results to traditional questions and methods in symbolic dynamics.
Symbolic dynamics originated as a tool for studying
the sequence of state transitions (through discrete time) 
in the evolution (or trajectory) of a point in a dynamic system.

The study of dynamical systems often concentrates on a (compact)
metrizable space $S$ with a continuous transition map $\tau:S\to S$.  
This leads to the notion of
``topological conjugacy'' of such objects $(S,\tau)$, as we will discuss below.
First we describe the well-known topological and metric structure on
the set of walks in any graph.

For $\omega\in N(X)$, let $U_n(\omega)$ denote 
the set of all walks in $N(X)$ which agree with $\omega$ for the
first $n$ steps. This set depends only on 
the path given by the first $n$ steps of $\omega$; more precisely, $U_n(\omega)=U(\alpha)$,
where $\alpha=s_n(\omega)$ and $U(\alpha)$ denotes the preimage of $\alpha$ 
under the source truncation $s_n:N(X)\to P_n(X)$.  
Note that $U(\alpha)$ is empty unless $\alpha$ is
is the source truncation of some walk.

The sets $U(\alpha)$ are the ``cylinder sets''
used to study Markov chains and dynamical systems, as in
Kemeny and Snell [1976], Douglas and Lind [1995], Kitchens [1998], etc.
Note that, for any $\omega\in U_{n'}(\omega')\cap U_{n''}(\omega'')$, we have
$U_n(\omega)\subseteq U_{n'}(\omega')\cap U_{n''}(\omega'')$
where $n={\rm min}(n',n'')$.
This shows that the collection of all unions of sets of the form $U_n(\omega)$ 
is closed under arbitrary unions and finite intersections,
and thus gives a topology on $N(X)$.
We may refer to $N(X)$ with this topology as the {\it walk space} for graph $X$.

There is also a nice distance function on $N(X)$, given as follows:
let $d(\omega,\omega)=0$;
if $\omega$ and $\nu$ are distinct walks in $X$, and let $d(\omega,\nu)=2^{-n}$,
where $n$ is the smallest natural number such that $s_n(\omega)\neq s_n(\nu)$.
For example, we always have $d(\omega,\nu)\leq 1$; but
$d(\omega,\nu)< 1$ if and only if $d(\omega,\nu)\leq 1/2$ if and only if $s_0(\omega)=s_0(\nu)$
($\omega$ and $\nu$ have the same source node).
To show that this gives a metric on $N(X)$, we merely check the metric axioms:
$0=d(\omega,\nu)$ iff $\omega=\nu$, $d(\omega,\nu)=d(\nu,\omega)$, 
and $d(\omega,\nu)\leq d(\omega,\mu)+d(\mu,\nu)$,
for all $\omega,\nu,\mu$.  

In fact, $d$ satisfies the stronger {\it ultrametric condition}, 
$d(\omega,\nu)\leq {\rm max}(d(\omega,\mu),d(\mu,\nu))$
for all $\omega,\nu,\mu$, as is easy to check.
So the above distance function makes $N(X)$ an {\it ultrametric space}.  Then
$U_n(\omega)=\{\nu\in N(X): d(\omega,\nu)<2^{-n}\}$ is the open ball of radius $2^{-n}$ around $\omega$,
and the walk space topology has as its open sets the arbitrary unions
of open balls for the ultrametric.

\smallskip\noindent
{\bf Proposition:}

1) $N(X)$ is a totally disconnected topological space. 

2) $N(X)$ is a complete metric space for the ultrametric structure.

3) If $X$ is a finite graph, then $N(X)$ is compact and separable.

4) If $X$ has finitely many arcs leaving each node, then $N(X)$ is locally compact.

\smallskip\noindent
{\bf Sketch Proofs:}  For part 1, one shows that any subset of $N(X)$ with more than one element is not connected;
more precisely, if $\omega'\ne \omega$ then $\omega'\notin U_n(\omega)$, 
and $U_n(\omega)$ is open and closed.
For part 2, one constructs the limit of any cauchy sequence of walks.
For part 3, since $N(X)$ is metrizable, it suffices
to show that every sequence has a convergent subsequence; this is easy to do.
Also, $N(X)$ is separable since the periodic walks give a
countable dense set in it.
For part 4, one uses the fact that if $X(x,*)$ is finite for every node $x$, then the set of
paths of given length leaving $x$ is finite;
it follows that $U(\alpha)$ is a compact subspace of $N(X)$ for every path $\alpha$ of positive length. QED

\smallskip 
For example, if $X$ is a dynamic graph then $N(X)$ is a discrete topological space, since
 if $\omega$ and $\nu$ are distinct walks in the dynamic graph $X$,
then  $s_0(\omega)\neq s_0(\nu)$ and so 
$d(\omega,\nu)=2^{-0}=1$.  

On the other hand, for any set $S$, if $X=\B(S)$ (the bouquet with $S$ as its set of loops) then
the topology on $N(X)$ is the product
topology on $S^\N$, where $S$ is given the discrete topology.  

We have the following general results for the walk space topology.
The shift map $\tau:N(X)\to N(X)$ is continuous, since $\tau:N(X)\to N(X)$ satisfies
$d(\omega\circ\tau,\nu\circ\tau)\leq 2\cdot d(\omega,\nu)$ for all walks $\omega$ and $\nu$ in $N(X)$.
Also, if $f:X\to Y$ is a graph morphism, then $N(f):N(X)\to N(Y)$ is continuous, 
since  $N(f)$, as given by $\omega\mapsto f\circ\omega$, is ``distance decreasing'': 
$d(f\circ\omega,f\circ\nu)\leq d(\omega,\nu)$ for all walks $\omega$ and $\nu$ in $X$.

\smallskip\noindent
{\bf Definition:}  A graph morphism $f:X\to Y$ is a {\it topological N-equivalence} if $N(f)$ is a
topological conjugacy.
Graphs $X$ and $Y$ are {\it topologically N-equivalent}  
(denoted  $X\sim_{tN} Y$)
when there exists an isomorphism of N-spaces $\phi:N(X)\to N(Y)$ which is a homeomorphism.
Then $N(X)$ and $N(Y)$ are said to be {\it topologically conjugate},
and $\phi$ is said to be a {\it topological conjugacy}. 

\smallskip\noindent
{\bf Proposition:} For any graph $X$, the graph morphism $s_{m,n}:\e^{n+m}X\to \e^nX$ is a 
topological N-equivalence for all $n,m\geq 0$.
In particular, $s_{n,0}:\e^nX\to X$ is a topological N-equivalence.  But $s_n:\W X\to \e^nX$
is not in general a topological N-equivalence.

\smallskip\noindent
{\bf Proof:}   We have already shown that $N(s_{n,m}):N(\e^{n+m}X)\to N(\e^mX)$ is an isomorphism of N-sets.
Consider $s_{n,0}:\e^n X\to X$. To show that $N(s_{n,0})$ is a homeomorphism,
we observe that if walks $\omega$ and $\nu$ in $X$
correspond to walks $\omega'$ and $\nu'$ in $\e^nX$, then
$d(\omega',\nu')=k\cdot d(\omega,\nu)$, where $k=2^{n}$. 
The first statement follows when we replace $X$ by $\e^nX$.
Taking $X=\B(\{ a,b\})$, the bouquet on two loops, shows that $N(s_0)$ is not a homeomorphism,
since the topological space $\{ a,b\}^{\N}$ is not discrete,
while $N(\W X)$ has the discrete topology for any graph $X$. QED 

\smallskip
It follows that any finite-level homotopy arrow between graphs gives a continuous N-set map,
since $N(\e^nX)\to X$ is a homeomorphism and $N(f):N(\e^nX)\to N(Y)$ is continuous.
So the equivalence of categories from Ho(Gph) to NSet actually carries the finite-level homotopy subcategory 
into a topologized category of N-sets.  In particular, we have the following.

\smallskip\noindent
{\bf Corollary:}  If graphs $X$ and $Y$ are finite-level homotopy-equivalent 
then they are topologically N-equivalent.

\smallskip
For finite graphs we have the following result, of the type
attributed to  Curtis, Lyndon, and Hedlund 
in Lind and Marcus [1995] (page 186); they use the terminology ``finite-type shift space''
for $N(X)$, and ``sliding block code'' for $\phi$.

\smallskip\noindent
{\bf Proposition:} 
Let $X$ be a  finite graph.  If $\phi:N(X)\to N(Y)$ is a continuous N-map, 
then there exists a natural number $n$ and a graph morphism $f:\e^n X\to Y$
such that $\phi\circ N(s_{n,0})=N(f)$. 

\smallskip\noindent
{\bf Proof:}  Since $X$ is finite, the space $N(X)$ is compact;
so the continuous function $\phi:N(X)\to N(Y)$ is uniformly continuous.
In particular, there exists a constant $n$ so that, for every $\omega\in N(X)$,
$$\phi(U_n(\omega))\subseteq U_{0}(\phi(\omega))\quad
{\rm and}\quad \phi(U_{n+1}(\omega))\subseteq U_{1}(\phi(\omega)).$$
We define $f:\e^n X\to Y$ on nodes by $\alpha\mapsto s_0(\phi(\omega))$
where $\alpha=s_n(\omega)$;  and on arcs by $\beta\mapsto s_1(\phi(\omega))$
where $\beta=s_{n+1}(\omega)$.  The definition on nodes is independent of choice of $\omega$ since 
$s_n(\nu)=\alpha$ implies  $\nu\in U_n(\omega)$, which implies that
$\phi(\nu)\in U_{0}(\phi(\omega))$ and $s_0(\phi(\nu))=s_0(\phi(\omega))$.
The definition on arcs is similarly  independent of choice.  QED

We are most interested here in applying the above ideas to the
study of N-equivalence of finite graphs.

\smallskip\noindent
{\bf Proposition:} If $X$ a finite graph, then any N-equivalence $f:X\to Y$ is a topological N-equivalence.

\smallskip\noindent
{\bf Proof:} Since $X$ is finite, $N(X)$ is compact; and $N(Y)$ is metrizable
and thus Hausdorf.  So $N(f)$ is a continuous bijection which
carries closed sets to closed sets. Thus $N(f)$ is a homeomorphism. QED

\smallskip\noindent
{\bf Proposition:} For finite graphs $X$ and $Y$:

1)  $X$ and $Y$ are topologically N-equivalent
if and only if there exists a finite graph $E$ with N-equivalences $f:E\to X$ and $g:E\to Y$.

2) $X$ and $Y$ are topologically N-equivalent
if there exist N-equivalences $X\to B$, $Y\to B$.

\smallskip\noindent
{\bf Proof:} 
For part 1, if $f:E\to X$ and $g:E\to Y$ are N-equivalences with $E$ finite,
then $N(f)$ and $N(g)$ are topological conjugacies, so that $N(X)$ and $N(Y)$ are topologically conjugate.
Conversely, if $\phi:N(X)\to N(Y)$ is a topologically conjugacy, then the continuous map of N-sets
$\phi$ comes from some graph morphism $g:\e^n X\to Y$ (since $X$ is finite).
Let $E=\e^n X$; so $s:E\to X$ is an N-equivalence, and thus $g:\e^n X\to Y$
must be an N-equivalence since $\phi:N(X)\to N(Y)$ is an isomorphism.
For part 2, assume that $X\to B$ and $Y\to B$ are N-equivalences.
Consider the fiber-product (pullback) $E=X\times_B Y$.  
Then $E$ is finite, since $X$ and $Y$ are finite, and we have isomorphisms of N-sets
$$N(E)=N(X)\times_{N(B)} N(Y)=N(X)=N(Y).$$
So by part 1 we see that $X$ and $Y$ are topologically N-equivalent.  QED


\smallskip\noindent
{\bf 6. Some necessary and sufficient conditions for topological N-equivalence.}

In this section we want to give some new and old conditions for N-equivalence and
topological N-equivalence.  
In particular, we will explore connections between symbolic dynamics
and the following special type of graph morphism.  

\smallskip\noindent
{\bf Definition:} A graph morphism $f:X\to Y$ is a {\it covering}
when $f:X(*,x)\to Y(*,f(x))$ is a bijection 
for every node $x$ in $X$. We say it is an {\it epic covering} 
when $f$ is also surjective on nodes (and thus on arcs).
Here $X(*,x)$ denotes the set of arcs in $X$ with target the node $x$, etc.

According to the historical sketch given in Boldi and Vigna [2002], this
basic concept has independently arisen many times in graph theory. 
Other names for covering  include {\it divisor}, 
{\it fibration}, {\it equitable partition}, etc.
Many of the natural graph morphisms in this paper are coverings.

\smallskip\noindent
{\bf Proposition:} For any graph $X$, 
the source truncations $s_0:\W X\to X$ and $s_{n,0}:\e^n X\to X$ are coverings.
Also,  $s_n:\W X \to \e^n X$ and $s_{m,n}\e^{n+m} X\to \e^n X$  are coverings, for all $m,n\geq 0$.

\smallskip\noindent
{\bf Proof:} Since $f$ is an N-equivalence, $\W(f):\W X\to \W Y$ is a graph isomorphism.
A node in $\W X$ is a walk $\omega\in N(X)$.  Let $x=s(\omega)$.
Each arc in $\W X(*,\omega)$ has the form $(a\omega,a,\omega)$ with $a\in X(*,x)$;
so $\W X(*,\omega)\to X(*,x)$ is a bijection.
A similar argument applies to $\e^{n} X$, etc.
The final statement follows by applying the first results to the graph $\e^m X$. QED

\smallskip\noindent
{\bf Proposition:} If $X$ is walkable and $f:X\to Y$ is an N-equivalence then $f$ is a covering.

\smallskip\noindent
{\bf Proof:} Since $X$ is walkable, for any node $x$ in $X$ there is some  walk $\omega$ with source $x$.
Considering $\omega$ as a node in $\W X$, we have bijections $s_X:\W X(*,\omega)\to X(*,x)$
and $s_Y:\W Y(*,f\omega)\to Y(*,fx)$, and $\W(f):\W X(*,\omega)\to \W Y(*,f\omega)$. 
Moreover, $f\circ s_X=s_Y\circ\W(f)$.
It follows that $f:X(*,x)\to Y(*,fx)$ must be a bijection. QED

\smallskip
Recall that a graph morphism $f:X\to Y$ is 
a {\it level $n$ homotopy equivalence} when
there exists a graph morphism $q:\e^n Y\to X$  which fills the diagram  
$$          \diagram[nohug]
          \e^n X & \rTo^s   &X   \cr
  \dTo^{\e^n(f)} &  \NE^{q} &  \dTo^{f}   \cr
          \e^n Y & \rTo^s   &Y
                                    \enddiagram$$

\smallskip\noindent
{\bf Proposition:} If $f:X\to Y$ is a level $n$ homotopy equivalence
and every node is the source of 
some path $\alpha$ of length $n$ in $X$, then
$f$ is a covering. 

\smallskip\noindent
{\bf Proof:} By hypothesis,
for any node $x$ in $X$ 
there is some path $\alpha$ of length $n$ with source $x$.
The graph morphism $q:\e^n Y\to X$ satisfies $q\circ \e^n(f)=s_X$ and $f\circ q=s_Y$. 
Considering $\alpha$ as a node in $\e^n X$ and $f\alpha$ as a node in $\e^n Y$,
we have bijections $s_X:\e^n X(*,\alpha)\to X(*,x)$
and $s_Y:\e^n Y(*,f\alpha)\to Y(*,fx)$.
Since $q\circ \e^n(f)=s_X$, we know that $q(f\alpha)=x$. 
Consider the function $q:\e^n Y(*,f\alpha)\to X(*,x)$. 
Since $f\circ q=s_Y:\e^n Y(*,f(\alpha))\to Y(*,f(x))$ is a bijection, 
we know that $f:X(*,x)\to Y(*,f(x))$ is surjective. 
Since $q\circ \e^n(f)=s_X:\e^n Y(*,f(\alpha))\to Y(*,f(x))$ is a bijection, 
we know that $f:X(*,x)\to Y(*,f(x))$ is injective.
Thus $f:X(*,x)\to Y(*,fx)$ is a bijection. QED

\smallskip
We will use the above to derive a necessary condition for N-equivalence.

\smallskip\noindent
{\bf Definition:}  
Consider the graph $T(X,x)$ given as follows, where $x$ is a node in graph $X$. 
The nodes in $T(X,x)$ are the finite paths in $X$ with target $x$
(note that $x$ is considered as a path of length $0$ in $X$); 
the arcs in $T(X,x)$ are the triples
$(a\alpha,a,\alpha)$ where $a\alpha$ is the concatenation of path $\alpha$
and arc $a$ in $X$; and $s(a\alpha,a,\alpha)=a\alpha$
and $t(a\alpha,a,\alpha)=\alpha$.
There is a natural graph morphism $s:T(X,x)\to X$ given by $\alpha\mapsto s(\alpha)$ and
$(a\alpha,a,\alpha)\mapsto a$. 

\smallskip
The arcs in $T(X,x)$ which have the node $\alpha$  as target are those of the form
$(a\alpha,a,\alpha)$ for $a\in X(*,s(\alpha))$; it follows that
the graph morphism  $s:T(X,x)\to X$ is a covering.
Moreover, the graph $T(X,x)$ is a rooted tree, which we may call the {\it tree at $x$}.
Here by a {\it rooted tree}, we mean a graph $T$ with node $r$ such that
 there is a unique path in $T$ from $x$ to $r$, for each each node $x$ in $T$.
Notice that in this paper we are directing rooted trees {\it toward} their roots;
we used the opposite convention in Bisson, Tsemo [2008] and [2009].

An induction argument shows that if $f:X\to Y$ is a covering 
then $T(X,x)\to T(Y,f(x))$ is a graph isomorphism for every
node $x$ in $X$.
It follows that if  $f$ is a covering and nodes $x$ and $x'$ have $f(x)=f(x')$, then 
 $T(X,x)$ and $T(X,x')$ are isomorphic graphs.

\smallskip\noindent
{\bf Definition:} A graph $B$ is {\it basal} when the only epic coverings $B\to B'$ are isomorphisms.
A {\it basing for $X$} is an epic covering $p:X\to B$ where $B$ is basal.

\smallskip
The next three propositions are modeled on the discussion in Boldi and Vigna [2002].  
in their terminology, a basing is  a ``minimal fibration''.
We give the proofs here in our language (and with some added details).
We will refer to the graphs $T(B,x)$, for nodes $x$ in $B$, as the trees of $B$.

\smallskip\noindent
{\bf Proposition:} If no two trees in $B$ are isomorphic then $B$ is basal.

\smallskip\noindent
{\bf Proof:} If an epic covering is an injection on nodes then it must be an isomorphism.  So
if $p:B\to B'$ is an epic covering which is not an isomorphism, then there must be
at least two distinct nodes $x_1$ and $x_2$ in $B$ with $p(x_1)=p(x_2)$.  But this would
say that $B$ has two trees which are isomorphic.  
QED

\smallskip\noindent
{\bf Proposition:} Any graph $X$ has a basing $p:X\to B$.

\smallskip\noindent
{\bf Proof:} We define an equivalence relation on the nodes of $X$ by 
saying that nodes are equivalent when they have isomorphic trees.
Then we choose $B_0\subseteq X_0$ such that each equivalence class contains exactly one element of $B_0$.
Let $p_0:X_0\to B_0$ assign to each node in $X$ the element of $B_0$ in its equivalence class.
Define $B_1\subseteq X_1$ to be the disjoint union $B_1=\sum_{b\in B_0} X(*,b)$.  If we 
identify $B_1$ with the set of ordered pairs $(b,a)$ having $b\in B_0$ and $a\in X(*,b)$,
then we may define $s,t:B_1\to B_0$ by $s(b,a)=p_0(s(a))$ and $t(b,a)=b$.
The epic graph morphism $p:X\to B$ is given by function $p_0$ on nodes and by function $p_1(a)=(t(a),a)$ on arcs.
To show that $p$ is a covering, we use the bijection between $X(*,x)$ and $B(*,p(x))$ given by
the isomorphism between $T(X,x)$ and $T(X,p(x))$.
To show that $B$ is basal, we use the fact that 
if nodes $b,b'$ in $B$ have isomorphic trees, then the corresponding trees $T(X,b)$ and $T(X,b')$ are
isomorphic, so that $b=b'$. QED

\smallskip\noindent
{\bf Proposition:} If $B$ is basal then no two trees in $B$ are isomorphic.

\smallskip\noindent
{\bf Proof:} If two trees in $B$ were isomorphic, then the above construction would give an epic
covering $p:B\to B'$ which identifies the two nodes. This would not be an isomorphism, contradicting
the definition of basal graph. QED

\smallskip
Similar reasoning shows that if $p:X\to B$ is a basing and $X$ has isomorphic trees at nodes $x$ and $x'$, 
then $p(x)=p(x')$. We will use this in the next proof.  We will also use the notation
$f\approx_0 g$ to indicate that two graph morphisms $f$ and $g$ agree on nodes.

\smallskip\noindent
{\bf Proposition:} If $p:X\to B$ is a basing and $f:X\to Y$ is an epic covering then there exists 
an epic covering $h:Y\to B$ which ``commutes on nodes'' in that $(h\circ f)\approx_0 p$.

\smallskip\noindent
{\bf Proof:} Given an epic covering $f:X\to Y$ and a basing $p:X\to B$, we want to define a graph morphism
$h:Y\to B$ such that, on the level of nodes, $p_0=h_0\circ f_0$.  
Choose any section  $\phi:Y_0\to X_0$ for the surjective function $f_0:X_0\to Y_0$, 
so that  $f(\phi(y))=y$ for each node $y\in Y_0$.
Define $h$  on nodes by $h_0(y)=p(\phi(y))$; note that we have $p(x)=h_0(y)$ for any
node $x$ with $f(x)=y$, since then $T(X,x)$ is isomorphic to  $T(X,\phi(y))$, and $p$ is a basing.
But $\phi$ also determines a section $\phi_1:Y_1\to X_1$ of the surjective function $f_1:X_1\to Y_1$, 
by inverting each of the bijections $f:X(*,\phi(y))\to Y(*,y)$.  Define $h$ on arcs by $h_1(a)=p_1(\phi_1(a))$.
Let us check that this defines  a graph morphism  $h:X\to B$.
Let $y=t(a)$ and $y'=s(a)$; then
$$t(h(a))=t(p(\phi(a)))=p(t(\phi(a)))=p(\phi(y))=h(y)=h(t(a))$$
$$s(h(a))=s(p(\phi(a)))=p(s(\phi(a)))=p(\phi(y'))=h(y')=h(s(a))$$
note that $p(s(\phi(a))=p(\phi(y'))$ since $f(s(\phi(a)))=y'=f(\phi(y'))$.
In fact, $h$ is an epic covering since $h$ is surjection on nodes, and 
$h:Y(*,y)\to B(*,h(y))$, for each $y\in Y_0$, is the composition of bijections
$\phi_1:Y(*,y)\to X(*,\phi(y))$ and $p_1:X(*,\phi(y))\to B(*,p(\phi(y)))$. QED

\smallskip\noindent
{\bf Corollary:} If $p:X\to B$ and $p':X\to B'$ are basings then $B$ and $B'$ are isomorphic graphs.
More precisely, there exists an isomorphism of graphs $h:B'\to B$ with
$(h\circ p') \approx_0 p$.

\smallskip\noindent
{\bf Proof:}  The previous proposition, applied to the epic covering $p':X\to B'$
and the basing $p:X\to B$, gives the existence of an epic covering
$h:B'\to B$, which must be an isomorphism, since $B'$ is basal.  QED

\smallskip
So, we may speak of ``the basal graph of $X$'', as this is well-defined up to isomorphism of graphs.
But here is a cautionary example.

\smallskip\noindent
{\bf Example:}  Let $B=B'$ be the basal graph having one node $x$ and arcs $b,c$ 
(the bouquet with two loops).
Let $X$ have nodes $x_0$ and $x_1$ with arcs $b',b'',c',c''$ 
where 
$b':x_0  \to x_1$,
$b'':x_1  \to x_0$,
$c':x_0  \to x_0$, and
$c'':x_1  \to x_1$. 
Consider the graph morphism  $p:X\to B$ which takes $b',b''$ to $b$ and $c',c''$ to $c$,
and consider the graph morphism $p':X\to B'$ 
which takes $b',c'$ to $b$ and $b'',c''$ to $c$.
Note that  $p:X\to B$ and $p':X\to B'$ are epic coverings, and are thus basings;
but there is no graph morphism $f:B'\to B$ making $p=p'\circ f$.
So here are two basings $f:X \to B$ and $f':X \to B'$ 
which are not ``isomorphic'' (as graph morphisms),
even though their codomain basal graphs are isomorphic.

\smallskip\noindent
{\bf Proposition:} If $X$ and $Y$ are walkable graphs which are N-equivalent,
then the basal graphs of $X$ and $Y$ are isomorphic.

\smallskip\noindent 
{\it Proof:} Let $p:X\to B$ and $p':Y\to B'$ be basings for $X$ and $Y$.
Consider the source truncations $s:\W X\to X$ and $s:\W Y\to Y$.
These are coverings which are epic since $X$ and $Y$ are walkable.  
Since the graphs $X$ and $Y$ are N-equivalent, 
there exists a graph isomorphism $f:\W X\to \W Y$.
Thus we have epic coverings $p\circ s:\W X\to B$ and $p'\circ s\circ f:\W X\to B'$,
which are basings since $B$ and $B'$ are basal.  Thus $B$ and $B$ are ismorphic,
by the previous corollary. QED

\smallskip\noindent
{\bf Example:}
Note that each cycle graph $C_n$ has a basing to the terminal graph $1$,
but they are not N-equivalent  unless $n=1$, since their zeta series are different.
This shows that the converse of the above proposition is not true.

\smallskip
So isomorphism of basal graphs 
is a necessary condition for two graphs to be N-equivalent.
The following example shows that the basal graph is a finer invariant than the zeta series,
in that it can distinguish between N-equivalent graphs which have the same zeta series.

\smallskip\noindent
{\bf Example:} 
We exhibit two finite graphs 
which have the same zeta series but non-isomorphic basal graphs.  Let $X$ 
be the graph with nodes $0,1,2,3,4$ and arcs $(0,i)$ and $(i,0)$ for $i=1,2,3,4$. Let $Y$ 
be the graph with nodes the integers mod $4$,
with arcs $(i,i+1)$ and $(i,i-1)$ for all $i$ mod $4$,
and with source and target given by $s(i,j)=i$ and $t(i,j)=j$.
The characteristic polynomial of $Y$ is $x^4-4x^2$ 
and the characteristic polynomial of $X$ is $x^5-4x^3$;
so $X$ and $Y$ have the same zeta series (see the discussion at the end of Bisson and Tsemo [2009]).
But $X$ has a basing to the graph $B$ with nodes
$x$ and $x'$ and with four arcs from $x$ to $x'$ and one arc from $x'$ to $x$;
while $Y$ has a basing to the graph $B'$ with one node and two loops.
Since $B$ does not have the same number of nodes as $B'$, it follows that $X$ and $Y$ are not
N-equivalent (so that $N(X)$ and $N(Y)$ are not isomorphic as N-sets).

\smallskip 
Two arcs $a$ and $a'$ in graph $Y$ are said to be {\rm parallel}
when $s(a)=s(a')$ and $t(a)=t(a')$.
A graph $Y$ is said to be {\rm separated} when it has no parallel arcs.
This terminology comes from Vigna [1997], where he discusses
some of the features of the full subcategory of separated graphs.
We need the following simple observation:  
if $Y$ is a separated graph, then graph morphisms $f,g:X\to Y$ are equal 
if and only if $f\approx_0 g$.

\smallskip\noindent
{\bf Proposition:} If $X$ and $B$ are finite and walkable, and $B$ is separated and basal, 
then $X$ and $B$ are topologically N-equivalent if and only if they are N-equivalent.

\smallskip\noindent
{\bf Proof:} Clearly topological N-equivalence implies N-equivalence.
Assume that $X$ and $B$ are N-equivalent graphs which are finite and walkable;
and assume also that $B$ is separated and basal.
So we have a graph isomorphism $f:\W X\to \W B$ and $s:\W B\to B$ is a basing
(since $B$ is walkable, $s$ is an epic covering).  So $s\circ f:\W X\to B$ is a basing.
Let $p':X\to B'$ be a basing.  Then $s:\W X\to X$ is an epic covering since $X$ is
walkable, and so $p'\circ s:\W X\to B'$ is a basing.  Since $s\circ f$ and $p'\circ s$
are both basings of $\W X$, it follows that there exists an isomorphism of graphs
$h:B'\to B$ such that $(h\circ p'\circ s)\approx_0 (s\circ f)$.  But $B$ is separated, so
we must have $h\circ p'\circ s = s\circ f$.  
Since $s:\W X\to X$ and $s\circ f:\W X\to B$ are N-equivalences,
the graph morphism $h\circ p':X\to B$ must be an N-equivalence. Since
$X$ and $B$ are finite graphs, $h\circ p'$ must be a topological N-equivalence. QED

\vskip .5 in
\bigskip\noindent
{\bf Appendix A: Direct proof of the N-model structure on Gph.}
Here we show directly that our
three classes $({\cal C}_N,{\cal W}_N,{\cal F}_N)$ of graph morphisms,
from section 2,
satisfy the axioms for a model structure on Gph.

\smallskip
Let ${\cal W}_N$ be the N-equivalences.
Clearly ${\cal W}_N$ has the 2/3 property.

\smallskip
Let $\underline{\cal C}_N={\cal C}_N\cap{\cal W}_N$; this is the class Iso of isomorphisms in Gph.
It is clear that $(\underline{\cal C}_N,{\cal F}_N)=(\hbox{Iso},\hbox{All})$ is
a weak factorization system. 

\smallskip
Let $\underline{\cal F}_N={\cal F}_N\cap{\cal W}_N$; this is the class of N-equivalences,
$\underline{\cal F}_N={\cal W}_N$.  We must show that 
$({\cal C}_N,\underline{\cal F}_N)=({}^\dagger{\cal W}_N,{\cal W}_N)$ 
is a weak factorization system.
Let $f:X\to Y$ be an arbitrary graph morphism.
We must show that $f$ factors as $f=g\circ h$ with
$h\in {\cal C}_N$ and $g\in {\cal W}_N$. Recall the two graph morphisms
${\bf i}:0\to \N$ and ${\bf j}:\N+\N\to \N$, which are easily seen to
be in ${\cal C}_N={}^\dagger{\cal W}_N$.
We will construct $h$ as a transfinite composition of pushouts of copies of
${\bf i}$, and ${\bf j}$, from which $h\in {\cal C}_N$ follows, by general principles.
We will give a complete description of the construction here, since it involves
a ``small object argument''
(for these ideas, see Section 2.1 in Hovey [1999], for instance).

\smallskip
First we produce a graph $X'$ and graph morphisms $f':X\to X'$ and $g':X'\to Y$,
with $f=g'\circ f'$, and with $f'\in {\cal C}_N$ and $N(g')$ a surjection.
We construct $f'$ as a pushout of copies of ${\bf i}$, as follows.
For any set $I$ we can form a graph morphism 
$$\sum_I {\bf i}:\sum_I 0\to\sum_I \N$$ 
Take $I=N(Y)$, which determines a unique graph morphism $k:(\sum_I \N)\to Y$.
Define $f':X\to X'$ by the pushout diagram
$$\diagram
          \sum_I 0 & \rTo^{\sum_I{\bf i}}& \sum_I \N  \cr
  		   \dTo^{} &                     &  \dTo^{k}   \cr
             X     & \rTo^{f'}           &X'
                                    \enddiagram$$
The graph morphisms $f:X\to Y$  and $k:(\sum_I \N)\to Y$ determine a unique graph morphism $g':X'\to Y$, 
with $g'\circ f'=f$.  We can see that $N(g')$ is a surjection, as follows. For any $\omega \in N(Y)$,
we have the inclusion $\omega':\N\to (\sum_I \N)$, which we may view as $\omega'\in N(\sum_I \N)$.
Let $\omega''$ be the image of $\omega'$ under $N(k):N(\sum_I \N)\to N(X')$; then
$N(g')$ takes $\omega''$ to $\omega$.

\smallskip
Next we factor $g'$ through the composition of a number of steps.
We essentially use that every object in Gph is ``small'', 
and use a ``small object argument'' (following Section 2.1 in Hovey [1999]).
In fact, we may need a transfinite sequence of steps, 
so we will index our steps by a well-ordered set, an {\it ordinal}.
Take each ordinal to be the set of all smaller ordinals 
(see Chapter II, Section 3 in Cohen [1966], for instance).
Then each ordinal $\alpha$ has a {\it successor}, defined as $\alpha+1=\alpha\cup\{\alpha\}$.

Let $\Lambda$ be an ordinal so large that there is no 
injective function $\Lambda\to X'_0\times X'_0$.
We also assume that $\Lambda$ is not the successor of any ordinal, so that
$\lambda\in\Lambda$ implies $\lambda+1\in\Lambda$. 
We view $\Lambda$ as a category
with an object for each element of $\Lambda$ and one morphism from $\lambda$ to $\lambda'$ when
$\lambda\leq\lambda'$, and we define a functor $X^\bullet:\Lambda\to {\rm Gph}$ equipped
with natural transformations $f^\bullet$ and $g^\bullet$.

We will define, for each $\lambda \in \Lambda$, 
graph morphisms $f^\lambda:X'\to X^\lambda$ and $g^\lambda: X^\lambda\to Y$ 
with $g'=g^\lambda\circ f^\lambda$, and 
with $f^\lambda$ epic graph morphism in ${\cal C}_N$ 
and with $N(g^\lambda)$ surjective.
We actually define $X^\lambda$ and compatible graph morphisms 
$f^\lambda:X\to X^\lambda$ and $g^\lambda:X^\lambda\to Y$
by transfinite induction, assuming that they are defined for all smaller ordinals.
The transfinite inductive definition goes as follows.

\medskip
For the minimal element $0\in \Lambda$, let $X^0=X'$ and $f^0={\rm id}$ and $g^0=g'$, 
so that $g'=g^0\circ f^0$.
 
\smallskip
Assume that we have defined $X^\lambda$ and $f^\lambda$ and 
$g^\lambda$ with $f^\lambda\circ g^\lambda=g'$, for every $\lambda<\lambda'$, 
for some $\lambda'\in\Lambda$.

\smallskip
For $\lambda'$ a limit ordinal (not the successor of any ordinal) 
we define $X^{\lambda'}={\rm colim}_{\lambda<\lambda'}X^\lambda$.  The graph morphism
$f^{\lambda'}:X'\to X^{\lambda'}$, 
the {\it transfinite composition} of epimorphisms in ${\cal C}_N$, 
is an epimorphism in ${\cal C}_N$.
The colimit also determines a unique graph morphism
$g^{\lambda'}:X^{\lambda'}\to Y$, with $g'=g^{\lambda'}\circ f^{\lambda'}$.

\smallskip
If $\lambda'=\lambda+1$ and $N(g^\lambda)$ is  a bijection then 
we define $X^{\lambda +1}=X^\lambda$ and
$f^{\lambda +1}=f^\lambda$ and $g^{\lambda +1}=g^\lambda$.

\smallskip
If $\lambda'=\lambda+1$ and $N(g^\lambda)$ is not a bijection, then we define 
$X^\lambda\to X^{\lambda+1}$ by pushout with copies of ${\bf j}$,
indexed by the set $J$ of all $(\omega',\omega'')$ such that $N(g^\lambda)$
carries $\omega'$ and $\omega''$ to the same walk in $N(Y)$.

\medskip
We are gluing together along $(\N+\N)\to\N$ in each summand of $\sum_J (\N+\N)\to X^\lambda$, 
to produce an epimorphism $f^{\lambda +1}:X^\lambda\to X^{\lambda+1}$,
and a unique graph morphism $g^{\lambda'}:X^{\lambda'}\to Y$
with $g'=g^{\lambda'}\circ f^{\lambda'}$.
More precisely, $f^{\lambda +1}:X^\lambda\to X^{\lambda+1}$ is the pushout of 
$$\sum_J (\N+\N)\to X^\lambda\quad{\rm and}\quad \sum_J(\N+ \N)\to \N.$$

\smallskip\noindent
Note that if $g^\lambda$ is an N-equivalence, then we will have $X^\lambda=X^{\lambda'}$ 
for all $\lambda'>\lambda$,
and we may say that the {\it $\Lambda$-sequence stabilizes at $\lambda$}.
Let us verify that our $\Lambda$-sequence stabilizes at some $\lambda\in \Lambda$,
so that $g'=g^\lambda\circ f^\lambda$; then $f=g\circ h$ 
with $h=f^\lambda\circ f'$ and $g=g^\lambda$
gives our desired factorization, with $h\in{\cal C}_N$ and $g\in {\cal W}_N$.

Each graph epimorphism $f^\lambda:X'\to X^\lambda$ determines an equivalence relation
$E^\lambda\subseteq X'_0\times X'_0$ on the nodes of $X'$.
So long as $g^\lambda$ is not an N-equivalence,
we have $E^\lambda\subset E^{\lambda+1}$, a strict
inclusion.  This shows that the $\Lambda$-sequence 
constructed above eventually stabilizes, since otherwise we could choose a 
$\Lambda$-parametrized family of elements $p^{\lambda}\in X'_0\times X'_0$
with $p^{\lambda+1}\in E^{\lambda+1}-E^\lambda$.  This would give an injective function
$\Lambda\to X'_0\times X'_0$, which is impossible by our assumption about the size of $\Lambda$.
QED

\bigskip

\centerline{\bf Bibliography}

\smallskip
\item{[2003]} C. Berger  and I. Moerdijk, Axiomatic homotopy theory for operads. 
Comment. Math. Helv. 78 (2003), no. 4, 805--831.

\smallskip
\item{[2008]} T. Bisson and A. Tsemo, A homotopical algebras of graphs
related to zeta series, Homology, Homotopy and its Applications, 10 (2008),
1-13. 

\smallskip
\item{[2009]} T. Bisson and A. Tsemo, Homotopy equivalence of isospectral graphs,
on the arXiv since June 2009.

\smallskip
\item{[2002]} P. Boldi and S. Vigna, Fibrations of graphs,
Discrete Math., 243 (2002), 21-66.

\smallskip
\item{[1966]} P. J. Cohen, {\it Set Theory and the Continuum Hypothesis}, W.A. Benjamin, NY, 1966.

\smallskip
\item{[1992]}  D. B. Epstein, {\it Word Processing in Groups}, AK Peters, 1992

\smallskip
\item{[1999]} M. Hovey, {\it Model Categories}, Amer. Math. Soc., Providence, 1999.

\smallskip
\item{[2007]} A. Joyal and M. Tierney, Quasi-categories vs Segal spaces, 277-326 in
{\it Categories in algebra, geometry and mathematical physics}, Contemp. Math.
431, Amer. Math. Soc., Providence, 2007. 

\smallskip
\item{[1976]} J.G. Kemeny and J.L. Snell and A.W. Knapp, {\it Denumerable Markov Chains}.
Springer-Verlag, New York Berlin Heidelberg, 1976.

\smallskip
\item{[1998]} B.P. Kitchens, {\it Symbolic Dynamics: One-sided, Two-sided and Countable State Markov Shifts}.
Universitext, Springer-Verlag, New York Berlin Heidelberg, 1998.

\smallskip
\item{[2000]} M. Kotani and T. Sunada, Zeta functions of finite graphs,
J. Math. Sci. Univ. Tokyo, 7 (2000), 7-25.

\smallskip
\item{[1989]} F.W. Lawvere, Qualitative distinctions between some toposes of
generalized graphs, 261-299 in {\it Categories in computer science and logic
(Boulder 1987)}, Contemp. Math. 92, Amer. Math. Soc., Providence, 1989.

\smallskip
\item{[1997]} F.W. Lawvere and S.H. Schanuel, {\it  Conceptual Mathematics: a first
introduction to categories}. Cambridge University Press, Cambridge, 1997. 

\smallskip
\item{[1995]} D. Lind and B. Marcus,  {\it An introduction to symbolic dynamics and
coding}. Cambridge University Press, Cambridge, 1995. 

\smallskip
\item{[1971]} S. Mac Lane,  {\it Categories for the working mathematician},
Graduate Texts in Mathematics, Vol. 5. Springer-Verlag, New York-Berlin, 1971 (Second
edition, 1998).

\smallskip
\item{[1994]} S. Mac Lane and I. Moerdijk, {\it Sheaves in Geometry and Logic: a
first introduction to topos theory}, Universitext, Springer-Verlag, New York
(1994).

\smallskip
\item{[1967]} D.G. Quillen, {\it Homotopical Algebra}, Lecture Notes in Mathematics no. 43, Springer-Verlag, Berlin, 1967.

\smallskip
\item{[2005]} I. Raeburn,  {\it Graph algebras},
Conference Board of the Mathematical Sciences, 
the American Mathematical Society, Providence, R.I.  2005.

\smallskip
\item{[2009]} J. Sakarovitch, {\it Elements of Automata Theory}, (tr. R. Thomas),
Cambridge University Press, New York,  2009.

\smallskip
\item{[1997]} S. Vigna, A guided tour in the topos of graphs. Technical Report 199-97.
Universit\`a di Milano. Dipartimento di Scienze dell'Informazione, 1997.

\bigskip\bigskip\noindent
Terrence  Bisson, \quad bisson@canisius.edu

Department of Mathematics and Statistics, Canisius College, 
2001 Main Street, Buffalo, NY 14216 USA

\bigskip\noindent
Aristide Tsemo, \quad tsemo58@yahoo.ca

102 Goodwood Park, Apt. 614, Toronto, Ontario M4C 2G8 Canada

\end